\documentclass[11pt]{article}

\usepackage{amsmath,amsfonts,amsthm,amssymb,graphicx,bm}
\usepackage{enumerate}

\usepackage{bbm}

\numberwithin{equation}{section}

\newcommand{\RR}{\mathbb{R}}
\newcommand{\Sph}{\mathbb{S}}
\newcommand{\eps}{\epsilon}

\newcommand{\sign}{\mathrm{sign}\,}
\newcommand{\ind}[1]{\mathbf{1}_{\{#1\}}\,}

\newcommand{\footnoteremember}[2]{

  \footnote{#2}
  \newcounter{#1}
  \setcounter{#1}{\value{footnote}}

} \newcommand{\footnoterecall}[1]{

  \footnotemark[\value{#1}]

} 

\title{Mathematical description of bacterial traveling pulses}

\author{
Nikolaos Bournaveas\protect\footnoteremember{1}{University of Edinburgh, School of Mathematics,
JCMB, King's Buildings, Edinburgh EH9 3JZ, UK.
Email: \texttt{n.bournaveas@ed.ac.uk}} , 
Axel Buguin\protect\footnoteremember{2}{Institut Curie, UMR CNRS 168 "Physico-Chime-Curie"
11  rue Pierre et Marie Curie,
F-75248 Paris cedex 05, France. Emails: \texttt{axel.buguin@curie.fr}, \texttt{jonathan.saragosti@curie.fr}, \texttt{pascal.silberzan@curie.fr}} , 
Vincent Calvez\protect\footnoteremember{3}{\'Ecole Normale Sup\'erieure de Lyon, UMR CNRS 5669 "Unit\'e de Math\'ematiques Pures et Appliqu\'ees", 
46 all\'ee d'Italie, F-69364 Lyon Cedex 07,
France. Email: \texttt{vincent.calvez@ens-lyon.fr}} , 
\\
Beno\^it Perthame\protect\footnoteremember{4}{Universit\'e Pierre et Marie Curie-Paris 6, UMR CNRS 7598 "Laboratoire Jacques-Louis Lions", BC187, 4  place Jussieu,  F-75252 Paris cedex 5, and Institut Universitaire de France. Email: \texttt{benoit.perthame@upmc.fr}} , 
Jonathan Saragosti\protect\footnoterecall{2} ,
Pascal Silberzan\protect\footnoterecall{2}
}

\date{}

\begin{document}

\maketitle

\begin{abstract}
The Keller-Segel system has been widely proposed as a model for bacterial waves driven by chemotactic processes. Current experiments on {\em E. coli} have shown precise structure of traveling pulses. We present here an alternative mathematical description of traveling pulses at a macroscopic scale. This modeling task is complemented with numerical simulations in accordance with the experimental observations.  
Our model is derived from an accurate kinetic description of the mesoscopic run-and-tumble process performed by bacteria. This model can account for recent experimental observations with {\em E. coli}. Qualitative agreements include the asymmetry of the pulse and transition in the collective behaviour (clustered motion versus dispersion). In addition we can capture quantitatively the main characteristics of the pulse such as the speed and the relative size of tails.
This work opens several experimental and theoretical perspectives. Coefficients at the macroscopic level are derived
from considerations at the cellular scale. For instance the stiffness of the signal integration process turns out to have a strong effect on collective motion. Furthermore the bottom-up scaling allows to perform preliminary mathematical analysis and write efficient numerical schemes. This model is intended as a predictive tool for the investigation of bacterial collective motion.\end{abstract}
\section{Introduction}

Since Adler's seminal paper \cite{Adler66}, several groups have reported the formation and the propagation of concentration waves in bacteria suspensions \cite{BudreneBerg,Park03}. Typically, a suspension of swimming bacteria such as {\em E. coli} self-concentrates in regions where the environment is slightly different such as the entry ports of the chamber (more exposed to oxygen) or regions of different temperatures. After their formation, these high concentration regions propagate along the channel, within the suspension. It is commonly admitted that chemotaxis (motion of cells directed by a chemical signal) is one of the key ingredients triggering the formation of these pulses. We refer to \cite{Tindall2008} for a complete review of experimental assays and mathematical approaches to model these issues and to 
\cite{BergBook} for all biological aspects of {\em E. coli}. 

\begin{figure}
\begin{center}
\includegraphics[width = .8\linewidth]{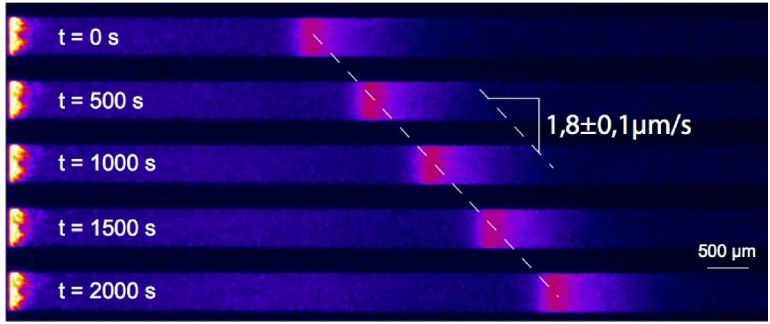}
\caption{\small Experimental evidence for pulses of {\em Escherichia coli} traveling across a channel. The propagation speed is constant and the shape of the pulse front is remarkably well conserved. Observe that the profile is clearly asymmetric, being stiffer at the back of the front (see also Fig. \ref{fig:front propagation}). Cell division may not play a crucial role regarding the short time scale.}
\label{fig:WaveChannel}
\end{center}
\end{figure}

Our goal  is to derive a macroscopic model for chemotactic pulses based on a mesoscopic underlying description (made of kinetic theory adapted to the specific run-and-tumble process that bacteria undergo \cite{Ajmb80,ODA}). We base our approach on recent experimental evidence for  traveling pulses (see Fig. \ref{fig:WaveChannel}). These traveling pulses  possess the following features which we are able to recover numerically: constant speed, constant amount of cells, short timescale (cell division being negligible), and strong asymmetry in the profile. 
  
We describe as usual the population of bacteria by its density $\rho(t,x)$ (at time $t>0$ and position $x\in \RR$). We restrict our attention to the one-dimensional case due to the specific geometry of the channels. The cell density follows a drift-diffusion equation, combining brownian diffusion together with  directed fluxes being the chemotactic contributions. This is coupled to reaction-diffusion equations driving the external chemical concentrations. In this paper we consider the influence of two chemical species, namely the chemoattractant signal $S(t,x)$, and the nutrient $N(t,x)$. Although this is a very general framework, it has been shown in close but different conditions that glycine can play the role of the chemoattractant \cite{Libchaber}. Similarly, glucose is presumed to be the nutrient. The exact nature of the chemical species has very little influence on our modeling process. In fact there is no need to know precisely the mechanisms of signal integration at this stage.
The model reads as follows:
\begin{equation}
\left\{
\begin{array}{rcl}
\partial_t \rho & = & D_\rho \Delta \rho - \nabla\cdot(\rho u_S + \rho u_N)\,, \medskip\\
\partial_t S & = & D_S \Delta S - \alpha S + \beta \rho\,, \medskip\\
\partial_t N & = & D_N \Delta N - \gamma \rho N\, .
\end{array}\right.
\label{eq:full model}
\end{equation}

The chemoattractant is assumed to be secreted by the bacteria (at a constant rate $ \beta$), and is naturally degraded at rate $\alpha$, whereas the nutrient is consumed at rate $\gamma$. Both chemical species diffuse with possibly different molecular diffusion coefficients. We assume a linear integration of the signal at the microscopic scale, resulting in a summation of two independent contributions for the directed part of the motion expressed by the fluxes $u_S$ and $u_N$. We expect that the flux $u_S$ will contribute to gather the cell density and create a pulse. The flux $u_N$ will be responsible for the motion of this pulse towards higher nutrient levels. Several systems of this type have been proposed and the upmost classical is the so-called Keller-Segel equation \cite{KellerSegel71,Murray} for which fluxes are proportional to the gradient of the chemical:
\[
u_S= \chi(S)\nabla S\, , \quad u_N= \chi(N)\nabla N\, .
\]
In the absence of nutrient, such systems enhance a positive feedback which counteracts dispersion of individuals and may eventually lead to aggregation. There is a large amount of literature dealing with this subtle mathematical phenomenon (see \cite{PerthameBook,HillenPainter} and references therein). Self-induced chemotaxis following the Keller-Segel model has been shown successful for modeling self-organization of various cell populations undergoing aggregation (slime mold amoebae, bacterial colony,\dots). 
In the absence of a chemoattractant $S$ being produced internally, this model can be used to describe  traveling pulses. However it is required that the chemosensitivity function $\chi(N)$ is singular at $N=0$ \cite{KellerSegel71}. 
Following the work of Nagai and Ikeda \cite{NagaiIkeda}, Horstmann and Stevens have constructed a class of such chemotaxis problems which admit traveling pulses solutions \cite{HorstmannStevens}, assuming the consumption of the (nutrient) signal together with a singular chemosensitivity. We also refer to \cite{Tindall2008} for a presentation of various contributions to this problem, and to \cite{LiWang} for recent developments concerning the stability of traveling waves in some parabolic-hyperbolic chemotaxis system.
%
%
In addition, the contribution of cell division to the dynamics of Keller-Segel systems (and specially traveling waves) has been considered by many authors (see \cite{Landman.TW,Mimura.TW,NadinPerthameRyzhik} and the references therein). However these constraints (including singular chemosensitivity or growth terms) seem unreasonable in view of the experimental setting we aim at describing.

An extension of the Keller-Segel model was also proposed in seminal paper by Brenner {\em et al.} \cite{BrennerLevitovBudrene} for the self-organization of {\em E. coli}. Production of the chemoattractant by the bacteria triggers consumption of an external field (namely the succinate). Their objective is to accurately describe aggregation of bacteria along rings or spots, 
as observed in earlier experiments by Budrene and Berg \cite{BudreneBerg}  that were performed over the surface of gels. One phase of the analysis consists in resolving a traveling ansatz for the motion of those bacterial rings. However the simple scenario they first adopt cannot resolve the propagation of traveling pulses. The authors give subsequently two possible directions of modeling: either observed traveling rings are transient, or they result from a switch in metabolism far behind the front. The experimental setting we are based on is quite different from Budrene and Berg's experiments (in particular regarding the dynamics): for the experiments discussed in the present paper, the bacteria swim in a liquid medium and not on agar plates. Therefore we will not follow \cite{BrennerLevitovBudrene}.

On the other hand Salman et al. \cite{Libchaber} consider a very similar experimental setting. However the model  they introduce to account for their observations is not expected to exhibit pulse waves (although the mathematical analysis would be more complex in its entire form than in \cite{HorstmannStevens}). Actually  Fig. 5 in \cite{Libchaber} is not compatible with a traveling pulse ansatz (because the pulse amplitude is increasing for the time of numerical experiments).

Traveling bands have also been reported in other cell species, and especially the slime mold {\em Dictyostelium discoideum} \cite{Weijer}.
Notice that the original model by Keller and Segel \cite{KellerSegel71} was indeed motivated by the observation of  traveling pulses in {\em Dictyostelium} population under starvation. This question has been developped more recently by H\"ofer et al. \cite{HoferSherrattMaini} using the Keller-Segel model, as well as Dolak and Schmeiser \cite{DolakSchmeiser} and Erban and Othmer \cite{ErbanOthmer07} using kinetic equations for chemotaxis. According to these models, the propagating pulse waves of chemoattractant (namely cAMP) are sustained by an excitable medium. The cells respond chemotactically to these waves by moving up to the gradient of cAMP. Great efforts have been successfully performed to resolve the "back-of-the-wave paradox": the polarized cells are supposed not to turn back when the  front passed away (this would result in a net motion outwards the pulsatile centers of cAMP). Although we are also focusing on the description of pulse waves, the medium is not expected to be excitable and the bacteria are not polarized.  Nevertheless, we will retain from these approaches the kinetic description originally due to Alt and coauthors \cite{Ajmb80,ODA}. This mathematical framework is well-suited for describing bacterial motion following a microscopic run-and-tumble process.

A new class of models for the collective motion of bacteria has emerged recently. It differs significantly from the classical Keller-Segel model. Rather than following intuitive rules (or first order approximations), the chemotactic fluxes are derived analytically from a mesoscopic description of the run-and-tumble dynamics at the individual level and possibly internal molecular pathways, see \cite{HillenOthmer,CMPS,FilbetLaurencotPerthame, ErbanOthmer04,CDMOSS,PerthameBook}. The scaling limit which links the macroscopic flux $u_S$ (or similarly $u_N$) to the kinetic description is now well understood since the pioneering works \cite{Ajmb80, ODA}. Here we propose to follow the analysis in \cite{DolakSchmeiser}, which is based on the temporal response of bacteria \cite{BrownBerg74,Segall86}, denoted by $\phi_\delta$ in Appendix \eqref{eq:kinflux}. Namely we write these fluxes as: 
\begin{equation}
\label{eq:genflux}
u_S =  \chi_S  J_\phi\left( -\epsilon\partial_t S,|\nabla S|\right)\dfrac{ \nabla S}{|\nabla S|} \, ,
\end{equation}
where $\epsilon$ is a (small) parameter issued from the microscopic description of motion. Namely $\epsilon$ is the ratio between the pulse speed and the speed of individual  cells
 (they differ by one order of magnitude at least according to experimental measurements).  
The function $J_\phi$ contains the microscopic features that stem from the precise response of a bacterium to a change in the environment (see the Appendix \eqref{eq:kinflux}). It mainly results from the so-called 'response function' at  the kinetic level that describe how a single bacterium responds to a change in the concentration of the chemoattractant $S$ in its surrounding environment. We give below analytical and numerical evidence that traveling pulses exist following such a modeling framework. 
We also investigate the characteristic features of those traveling pulses at the light of experimental observations.  

The experiments presented in the present paper will be described in more details in a subsequent paper. Briefly, in a setup placed under a low magnification fluorescence microscope maintained at $30\,^{\circ}\mathrm{C}$, we fill polymer microchannels -- section: ca. $500 \mu m\times 100 \mu m$, length ca. $1 cm$ -- with a suspension of fluorescent {\em E. coli} bacteria -- strain RP437 considered wild-type for motility and chemotaxis, transformed with a pZE1R-gfp plasmid allowing quantitative measurement of bacteria concentration inside the channel. We concentrate the cells at the extremity of the channel and monitor the progression of the subsequent concentration wave along the channel. In particular, we dynamically extract the shape of the front and its velocity (see Fig. \ref{fig:WaveChannel}).

Coupling the model \eqref{eq:full model} with the formula \eqref{eq:genflux} results in a  parabolic type partial differential equation for the bacterial density $\rho$, such as in the Keller-Segel system. It significantly differs from it however, as it derives in our case from a kinetic description of motion. Especially  the flux $u_S$ is uniformly bounded, whereas the chemotactic flux in Keller-Segel model generally becomes unbounded when aggregative instability occurs, which is a strong obstacle to the existence of traveling pulses.



\section{Traveling pulses under competing fluxes}
\label{sec:competing fluxes}

\subsection{Stiff response function: pulse wave analytical solutions}

\label{sec:analytical}

It is usually impossible to compute explicitely  traveling pulse solutions for general systems such as \eqref{eq:full model}. To obtain qualitative properties is also a difficult problem: we refer to \cite{NagaiIkeda, HorstmannStevens, NadinPerthameRyzhik} for examples of rigorous results. Here, we are able to handle analytical computations in the limiting case of a stiff signal response function, when the fluxes  are given  by the expression (see the Appendix \eqref{eq:kinflux}):
\begin{equation}
\label{eq:singfluxs}
u_S = \chi_S  \left(1 -  \left(\dfrac{ \epsilon\partial_t S}{ \partial_x S }\right)^2  \right)_+ \sign(\partial_x S)\, ,
\end{equation}
\begin{equation}
\label{eq:singfluxn}
u_N = \chi_N  \left(1 -  \left( \dfrac{\epsilon\partial_t N}{ \partial_x N }\right)^2  \right)_+ \sign(\partial_x N)\, .
\end{equation}
In other words, a specific expression for $J_\phi$ in \eqref{eq:genflux} is considered in this section. It eventually reduces to $u_S = \chi_S \sign(\partial_x S)$ as $\epsilon\to 0$.

We seek  traveling pulses, in other words particular solutions of the form $\rho(t,x) = \tilde\rho(x-\sigma t)$, $S(t,x) = \tilde S(x-\sigma t)$, $N(t,x) = \tilde N(x-\sigma t)$ where $\sigma$ denotes the speed of the wave. This reduces (\ref{eq:full model})   to a new system with a single variable $z=x-\sigma t$:
\begin{equation}
\label{eq:tw}
\left\{
\begin{array}{rcl} 
- \sigma \rho'(z) &=& D_\rho \rho''(z) - \left( \rho(z) u_S(z) + \rho(z) u_N(z) \right)'\, ,\medskip\\ 
- \sigma S'(z) &=& D_S S''(z) - \alpha S (z) + \beta \rho(z) \, ,  \medskip\\
- \sigma N'(z) &=& D_N N''(z) - \gamma \rho(z) N(z)\, .
\end{array}\right.
\end{equation}
We prescribe  the following conditions at infinity:
\begin{equation}
\label{eq:twcl}
\rho(\pm \infty )=0\,,  \qquad  S(\pm \infty )=0\,,  \qquad N(\pm \infty )=N_\pm \,.
\end{equation}


We impose $\sigma >0$ without loss of generality. This means that the fresh nutrient is located on the right side, and thus we look for an increasing nutrient concentration $N'(z)>0$.
We expect that the chemoattractant profile exhibits a maximum coinciding with the cell density peak (say at $z = 0$), and we look for a solution where  $S'(z)$ changes sign only once at $z = 0$. 
Then, the fluxes \eqref{eq:singfluxs}-\eqref{eq:singfluxn} express under the traveling wave ansatz as:
\begin{align*}
& u_S(z) =  - \chi_S  \left(1 -  \left(\epsilon{\sigma}{}\right)^2  \right)_+ \sign( z)\, ,\\
& u_N(z) =   \chi_N \left(1 -  \left(\epsilon{\sigma}{}\right)^2  \right)_+ \, .
\end{align*}

\begin{figure}
\begin{center}
\includegraphics[width = .8\linewidth, height = .4\textheight ]{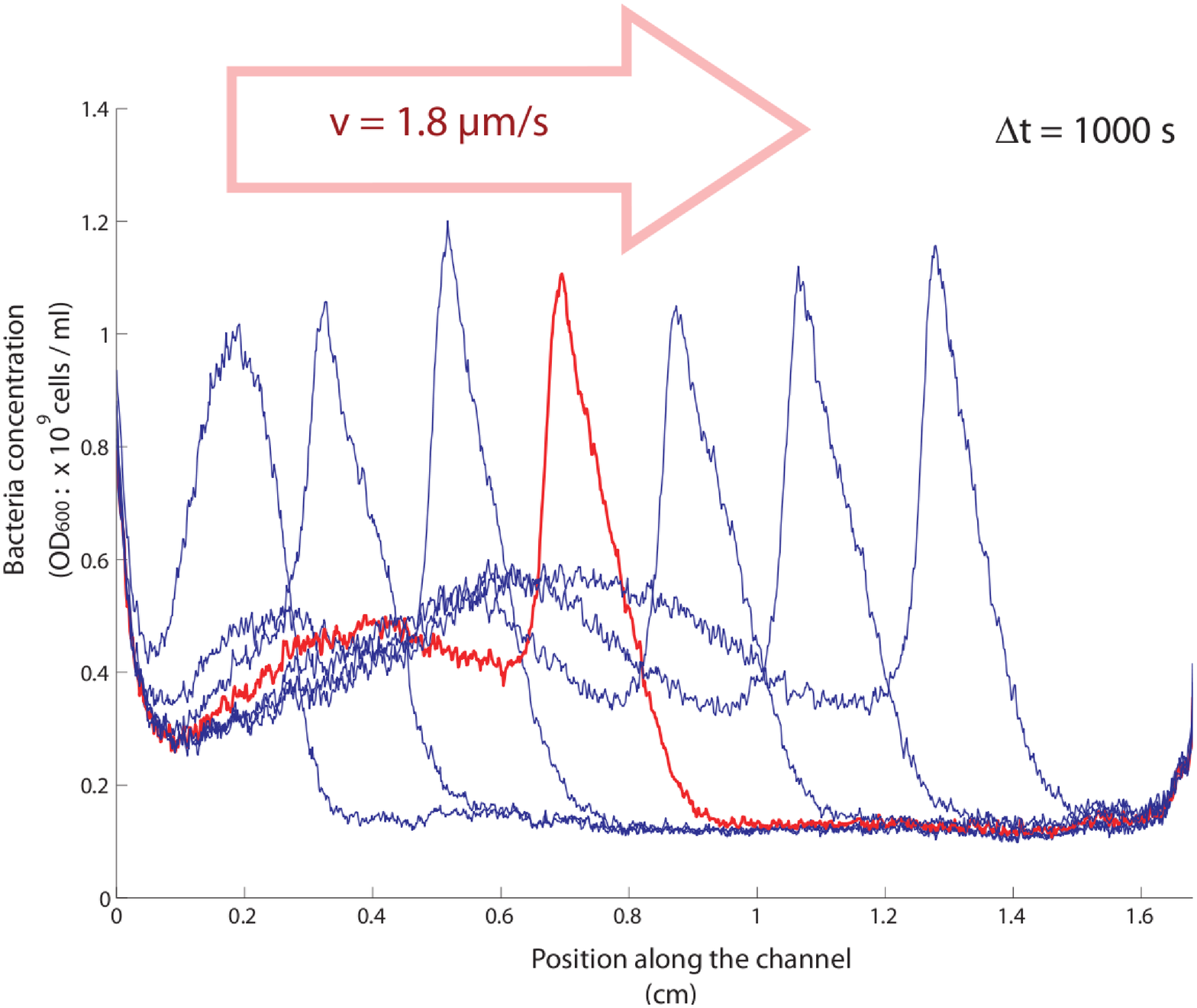}\\
\includegraphics[width = .6\linewidth]{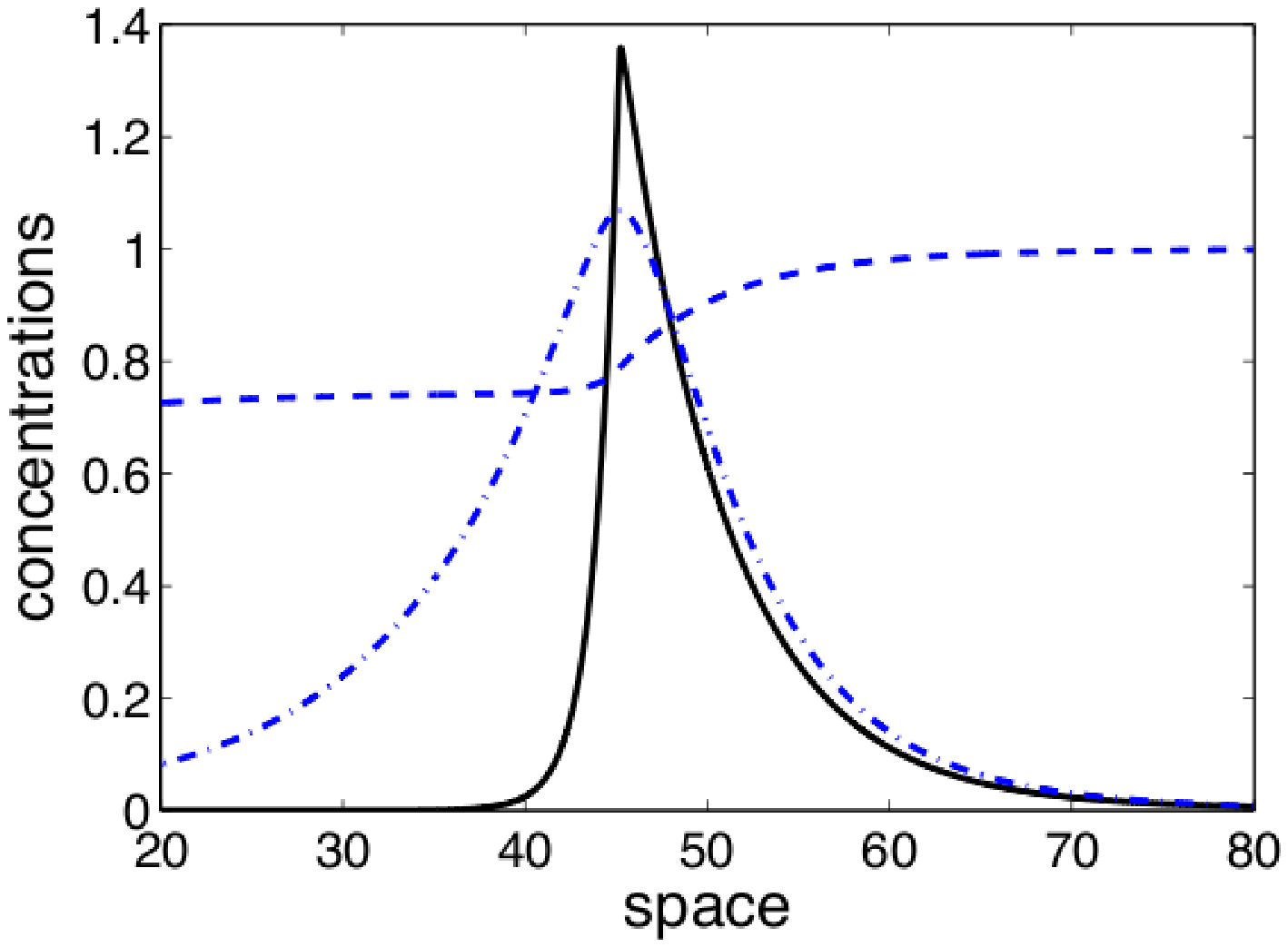}
\caption{\small (Top) Experimental evidence of a traveling pulse: time snapshots of the full experiments described in Fig. \ref{fig:WaveChannel}. The density profile is clearly asymmetric and preserved along the time course of the experiment. Main contribution to growth takes place at the back of the pulse. This suggests that nutrients are not totally consumed by the pulse. The number of bacteria in the pulse is approximately constant during the pulse course. (Bottom) A generic density profile obtained with the model \eqref{eq:full model} (see also Fig. \ref{fig:nolimited.sharp}). The bacteria (plain line) are attracted by a nutrient (dashed line) which is consumed. They secrete in addition their own chemical signal (dash-dotted line) which is concentrated at the peak location.}
\label{fig:front propagation}
\end{center}
\end{figure}

Integrating once the cell density equation we obtain,
\[
D_\rho \rho'(z) = \rho(z)\left(u_S(z) + u_N(z) - \sigma\right)\, .
\]
The flux $u_S$ takes two values (with a jump at $z=0$), whereas the flux $u_N$ is constant. Therefore the cell density is a combination of two exponential distributions:
\[
\rho(z) = \left\{
\begin{array}{lll}
\rho_0\exp\left(\lambda^- z\right)\, ,\,
& \lambda^- = \dfrac{ -\sigma +  (\chi_S+\chi_N)\left( 1 - \left(\epsilon{\sigma}\right)^2 \right) }{D_\rho}>0 \, ,\, & \mbox{if } z<0\, ,\\
\rho_0\exp\left(\lambda^+ z\right)\, ,\, 
& \lambda^+ =  \dfrac{-\sigma +  (-\chi_S+\chi_N)\left( 1 - \left(\epsilon{\sigma}\right)^2 \right)}{D_\rho}<0 \, ,\, & \mbox{if } z>0\, .
\end{array}
\right.
\]
This combination of two exponentials perfectly match with the numerical simulations (Fig. \ref{fig:front propagation}).

\subsection{Formula for the traveling speed} 
To close the analysis it remains to derive the two unknowns: the maximum cell density $\rho_0$ and the speed $\sigma$, given the mass and the constraint that $\partial_z S$ vanishes at $z = 0$ (because $S(z)$ reaches a maximum).

On the one hand, the total mass of bacteria is given by $M = \rho_0(1/\lambda^- +1/|\lambda^+|)$. On the other hand the chemotactic field is given by $S(z) = (K*\beta\rho)(z)$, where the fudamental solution of the equation for $S(z)$ is
\[ K(z) = a_1 \exp(-a_2|z| - a_3z)\, , \quad a_1 = \dfrac1{2 a_2 D_S}\,, \quad a_2  = \sqrt{a_3^2 + \dfrac{\alpha}{D_S}}\, , \quad a_3 = \frac{\sigma}{2D_S}\, . \]
To match the transition in monotonicity condition, the chemical signal should satisfy
$ S'(0) =  0$, that is $(K'*\beta\rho)(0)=0$, which leads to
\begin{eqnarray*}
0 &=& \rho_0 \int_{-\infty}^0 a_1(a_2+a_3)\exp( a_2z + a_3z)\exp(\lambda^- z)\, dz \\ &&\qquad + \rho_0 \int_0^{\infty} a_1(-a_2 + a_3)\exp(-a_2z+a_3z)\exp(\lambda^+ z)\, dz \\
0 & = & a_1  \left( \dfrac{a_2+a_3}{a_2+a_3+\lambda^-} - \dfrac{-a_2+a_3}{-a_2+a_3+\lambda^+} \right)\, .  
\end{eqnarray*}
This leads to the following equation that we shall invert to obtain the front speed:
\begin{eqnarray} \dfrac{\lambda^-}{|\lambda^+|} &=& \dfrac{a_2+a_3}{a_2-a_3}\, ,\nonumber\\
\dfrac{-\sigma +  (\chi_S+\chi_N)\left( 1 - \left({\epsilon\sigma} \right)^2 \right)}{\sigma +  (\chi_S-\chi_N)\left( 1 - \left({\epsilon\sigma}\right)^2 \right)} &=& 
\dfrac{\sqrt{4D_S\alpha+\sigma^2}+ \sigma}{\sqrt{4D_S\alpha+\sigma^2}-\sigma}\, .\label{eq:asymmetric factor} 
 \end{eqnarray}
From this relation we infer: 
\begin{equation}
\label{eq:front speed} 
\chi_N- \dfrac{\sigma}{1-\left({\epsilon\sigma} \right)^2}= \chi_S \dfrac{\sigma}{\sqrt{4D_S \alpha +\sigma^2}} \, .
\end{equation}
We deduce from monotonicity  arguments that there is a unique possible traveling speed $\sigma \in (0, \epsilon^{-1})$.


According to \eqref{eq:front speed} the expected pulse speed does not depend upon the total number of cells when the response function is stiff. This can be related to a recent work by Mittal {\em et al.} \cite{Mittal03} where the authors observe experimentally such a fact in a different context (see Section \ref{sec:cluster} below). In the case of a smooth tumbling kernel in \eqref{eq:kinflux}, our model would predict a dependency of the speed upon the quantity of cells. But this analysis suggests that the number of cells is presumably not a sensitive biophysical parameter. 

Observe that the speed does not depend on the bacterial diffusion coefficient either. Therefore we expect to get the same formula if we follow the hyperbolic approach of \cite{DolakSchmeiser} in order to derive  a macroscopic model. Indeed the main difference at the macroscopic level lies in the diffusion coefficient which is very small in the hyperbolic system. Nevertheless, the density distribution would be very different, being much more confined in the hyperbolic system. Furthermore, scaling back the system to its original variables, we would obtain a pulse speed being comparable to the individual speed of bacteria (see Appendix). This is clearly not the case.

\subsection{Cluster formation}
\label{sec:cluster}

\subsubsection*{Stiff response function.}

Mittal {\em et al.} have presented remarkable experiments where bacteria {\em E. coli} self-organize in coherent aggregated structures due to chemotaxis \cite{Mittal03}. The cluster diameters are shown essentially not to depend on the quantity of cells being trapped. This experimental observation can be recovered from direct numerical simulations of random walks \cite{Inoue07}.

We can recover this feature in our context using a model similar to \eqref{eq:full model} derived from a kinetic description.  Following Section \ref{sec:analytical} we  compute the solutions of \eqref{eq:tw} in the absence of nutrient (assuming again a stiff response function). Observe that stationary solutions correspond here to zero-speed traveling pulses. The problem is reduced to finding solutions of the following system: 
\begin{equation}
\label{eq:cluster}
\left\{
\begin{array}{l} 
- D_\rho \rho'(x) + \rho(x) u_S(x)  = 0 \, , \quad  u_S(x) =  \chi_S\sign(S'(x)) \, , \medskip\\ 
- D_S S''(x) + \alpha S (x) = \beta \rho(x) \, . 
\end{array}\right.
\end{equation}
We assume again that $\sign(S'(x)) = -\sign(x)$. This simply leads to,
\[
\rho(x)=\rho_0\exp(-\lambda |x|)\, , \quad\mbox{where}\quad  \lambda=\frac{\chi_S}{D_\rho}\, .
\]
This is compatible with the postulate that $S(x)$ changes sign only once, at $x=0$ (the source $\beta \rho(x)$ being even). The typical size of the clusters is of the order $\lambda^{-1}$, which does not depend on the total number of cells. This is in good agreement with experiments exhibited in \cite{Mittal03}. The fact that
we can recover them from numerical simulations indicates that these stationary states are expected to be stable.

\subsubsection*{General response function.}

Cluster formation provides a good framework for investigating the situation where we relax the stiffness assumption of the response function $\phi_\delta$. Below $\phi_\delta$ is characterized by the stiffness parameter $\delta$ through $\phi_\delta'(0) = -3/\delta$ (see Appendix). Consider the caricatural model (in nondimensional form):
\begin{equation}
\left\{\begin{array}{l}
\partial_t \rho =  \partial_{xx}^2 \rho - \partial_x\left( \rho u_S \right)\, , \displaystyle\quad  u_S = -\dfrac12 \int_{v\in (-1,1))} v \phi_\delta(\epsilon \partial_t S + v \partial_x S)\, dv\, , \medskip\\
-\partial_{xx}^2 S + \alpha S = \rho \, . 
\end{array}\right.
\end{equation}
We rewrite $\alpha = l^{-2}, $where $l$ denotes the range of action of the chemical signal. 
We investigate the linear stability of the constant stationary state $(\overline{\rho},\overline{S}) = (\langle\rho\rangle, \alpha^{-1}\langle\rho\rangle)$ where $\langle\rho\rangle$ is the meanvalue over the domain $[0,L]$. We introduce the deviation to the stationary state: $ n = \rho - \overline{\rho}$, $c = S - \overline{S}$. Then the linearized system writes close to $(\overline{\rho},\overline{S})$:
\begin{equation}
\left\{\begin{array}{l}
\partial_t n =  \partial_{xx}^2 n - \partial_x\left(\langle\rho\rangle  \widetilde{u_S} \right)\, , \quad \widetilde{u_S}  = \displaystyle -\dfrac12 \int_{v\in (-1,1)} v^2 \phi_\delta'(0)  \partial_x c \, dv   = \displaystyle \dfrac{ \partial_x c  }\delta\, , \medskip\\
-\partial_{xx}^2 c + \alpha c = n \, . 
\end{array}\right.  
\end{equation}
The associated eigenvalue problem reduces to the following dispersion relation for $\xi = 2\pi k /L$,
\[\lambda(k) = - \xi^2 + \dfrac{M}{\delta L} \dfrac{\xi^2}{\alpha + \xi^2}\, .\]
Due to the conservation of mass, we shall only consider $k\geq 1$. The eigenvalue becomes positive if there exists $k$ such that
\[ \dfrac M{\delta L} > \alpha + \xi^2\, , \quad \mbox{or equivalently}\quad 
\dfrac{M   }{2\pi\delta} > \dfrac{L}{2\pi l^2} + k^2 \dfrac{2\pi   }{L}\, .
\]
Therefore the constant solution is linearly stable if the following condition is fulfilled: \[ \dfrac{M L }{2\pi\delta} < \dfrac{1}{2\pi }\left(\dfrac L l\right)^2 +  2\pi \, . \]

The picture is not complete as we have not investigated the stability of the non-trivial steady-state. However this indicates that the mass and the stiffness parameter $\delta$ play important roles regarding cluster formation.



\section{Numerical insights}
\label{sec:num}

\begin{figure}
\begin{center}
\includegraphics[width = .6\linewidth]{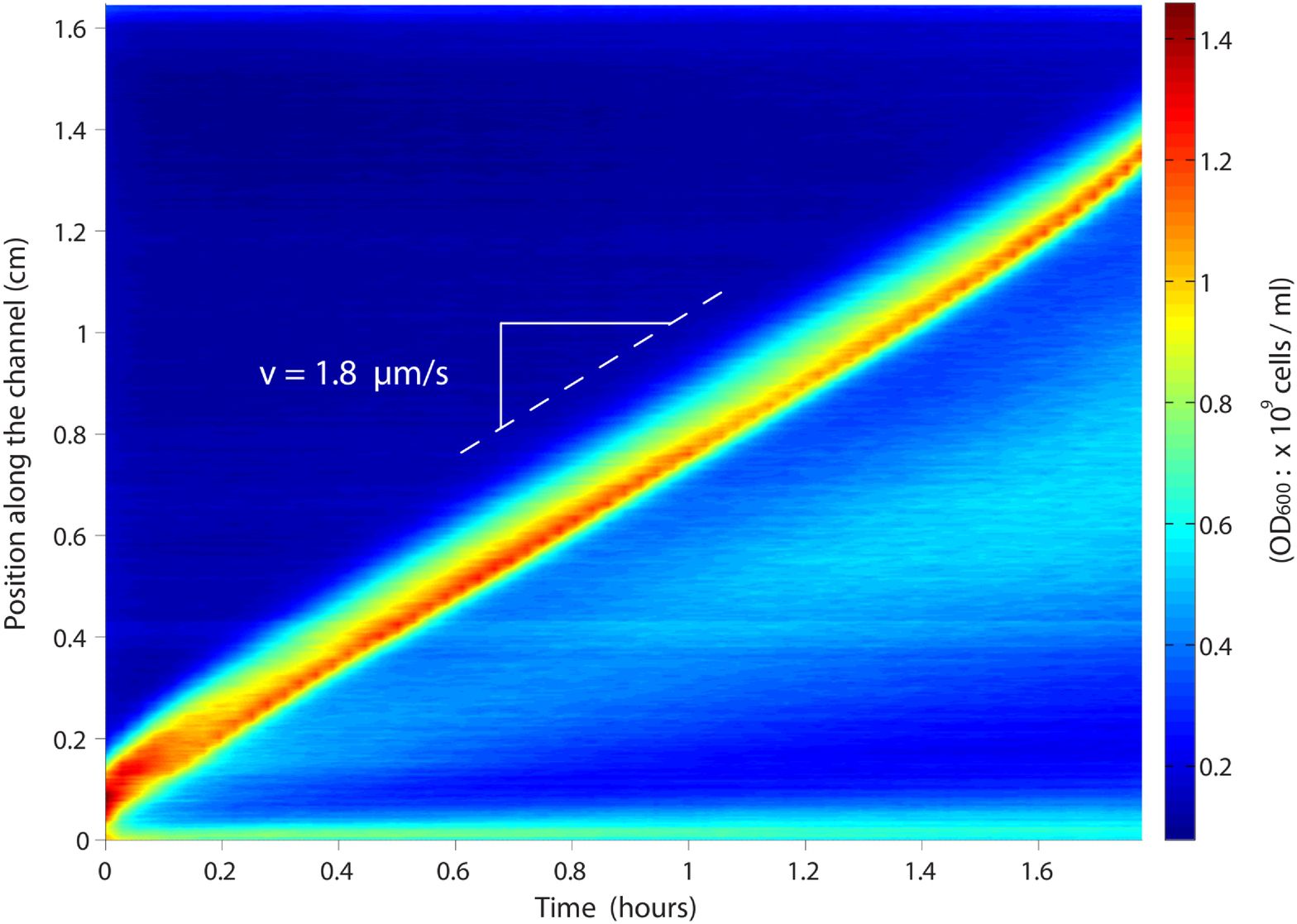}\\ 
\includegraphics[ width = \linewidth ]{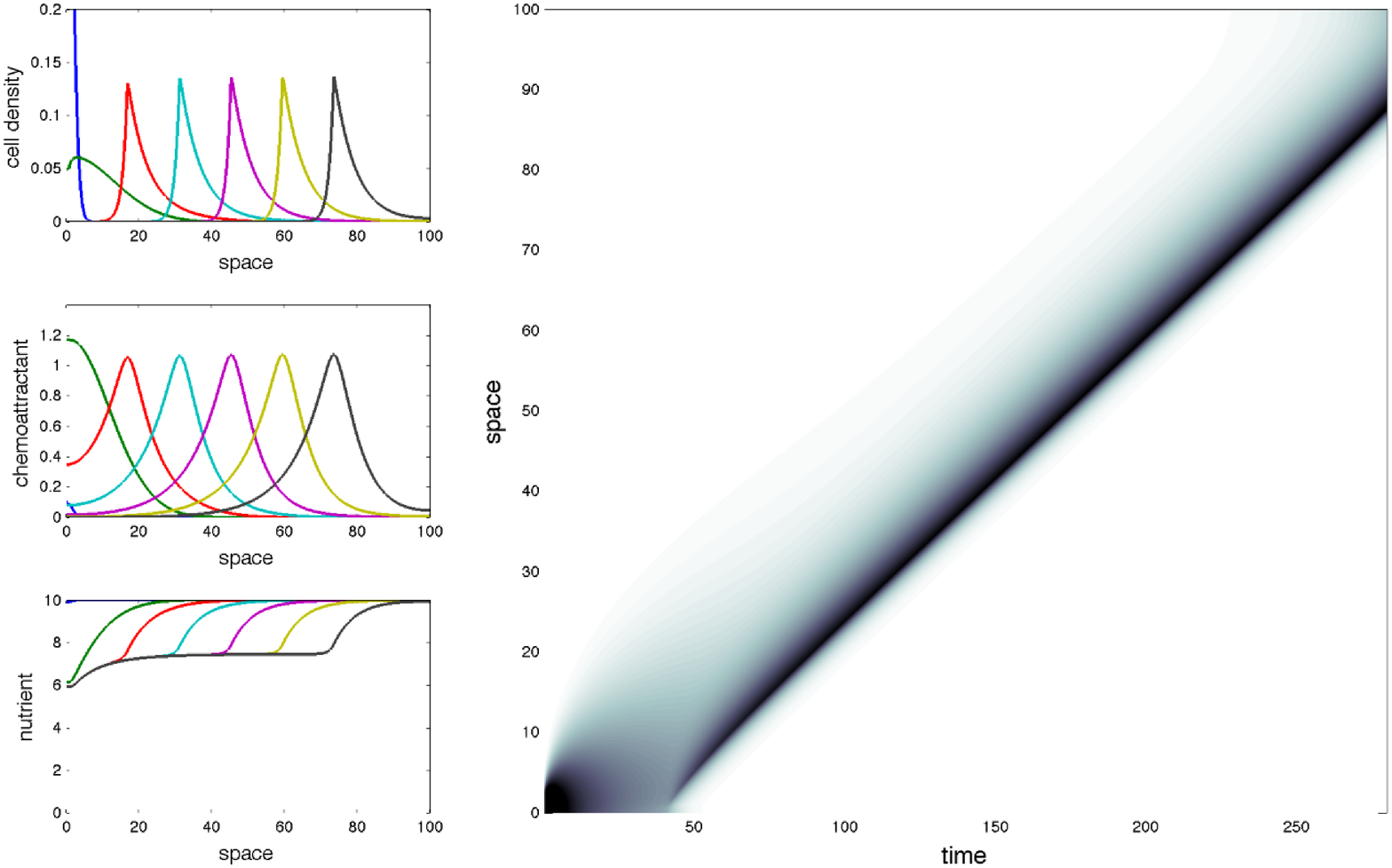}\caption{\small (Top) Experimental results under abundant nutrient conditions: M9 minimal medium supplemented with 4\% glucose and 1\% casamino acids (both ten times more concentrated than in the case of Fig. \ref{fig:limited nutrient}). (Bottom) Numerical simulations of system \eqref{eq:full model} in the case of unlimited nutrient, and a stiff response function $\phi_\delta$. We observe the propagation of a traveling pulse with constant speed and asymmetric profile. The value of the speed and the features of the profile (combination of two exponential tails) perfectly match the traveling wave ansatz analysis of Section \ref{sec:competing fluxes}. We use a semi-implicit upwind finite-difference scheme performing a half-point discretization (in space) of the time derivative of chemical species when computing the approximate fluxes. Specific parameters are: $\delta = 10^{-3}$ and $N_0 = 10$ (see Section \ref{sec:num} for the other parameters). 
}
\label{fig:nolimited.sharp}
\end{center}
\end{figure}

We complete the theoretical analysis with some numerical simulations exhibiting propagation of pulses (or not)
in regimes where analytical solutions are not available (see Fig. \ref{fig:nolimited.sharp}).


\begin{figure}
\begin{center}
\includegraphics[width = .6\linewidth]{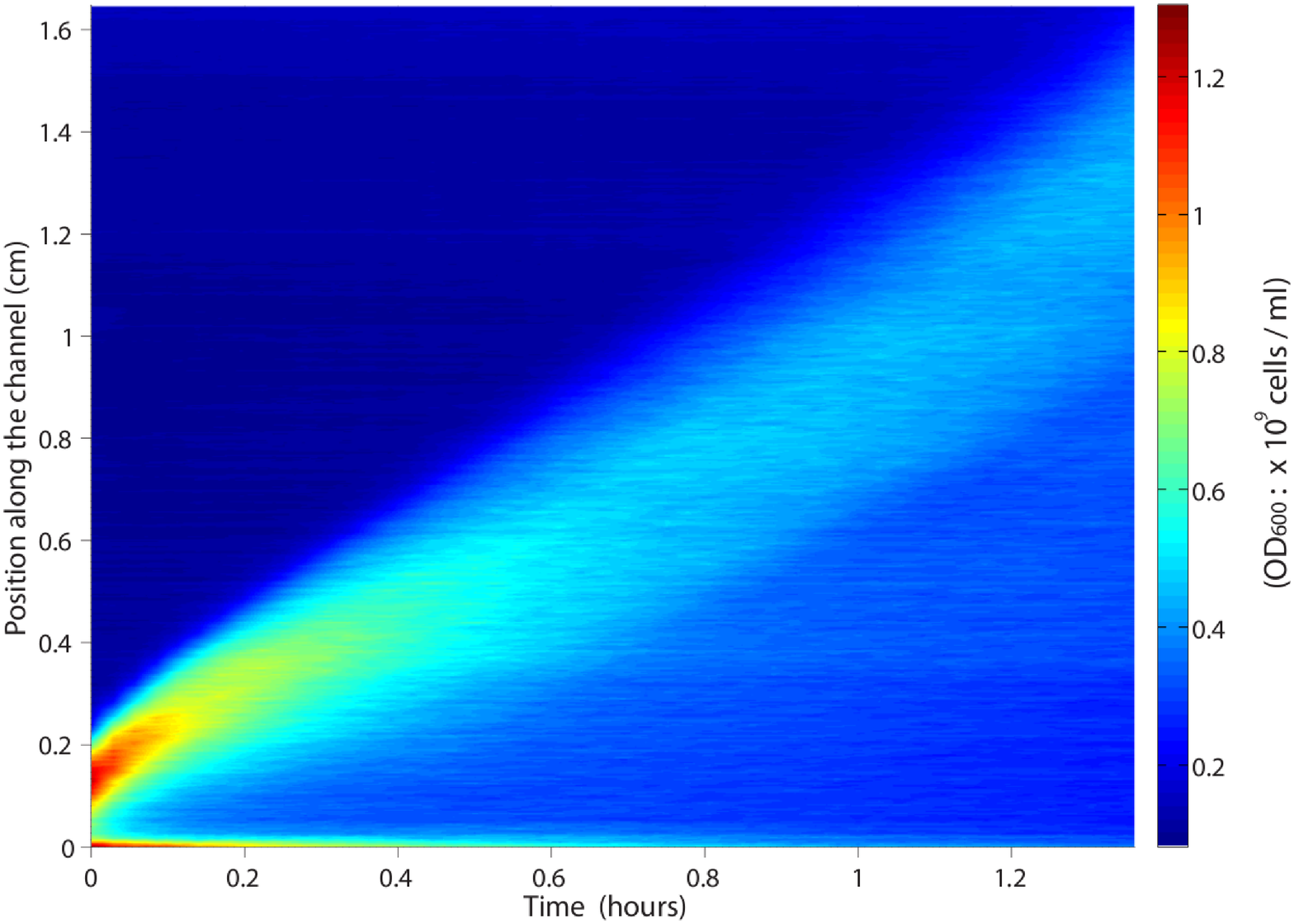}\\
\includegraphics[width= \linewidth]{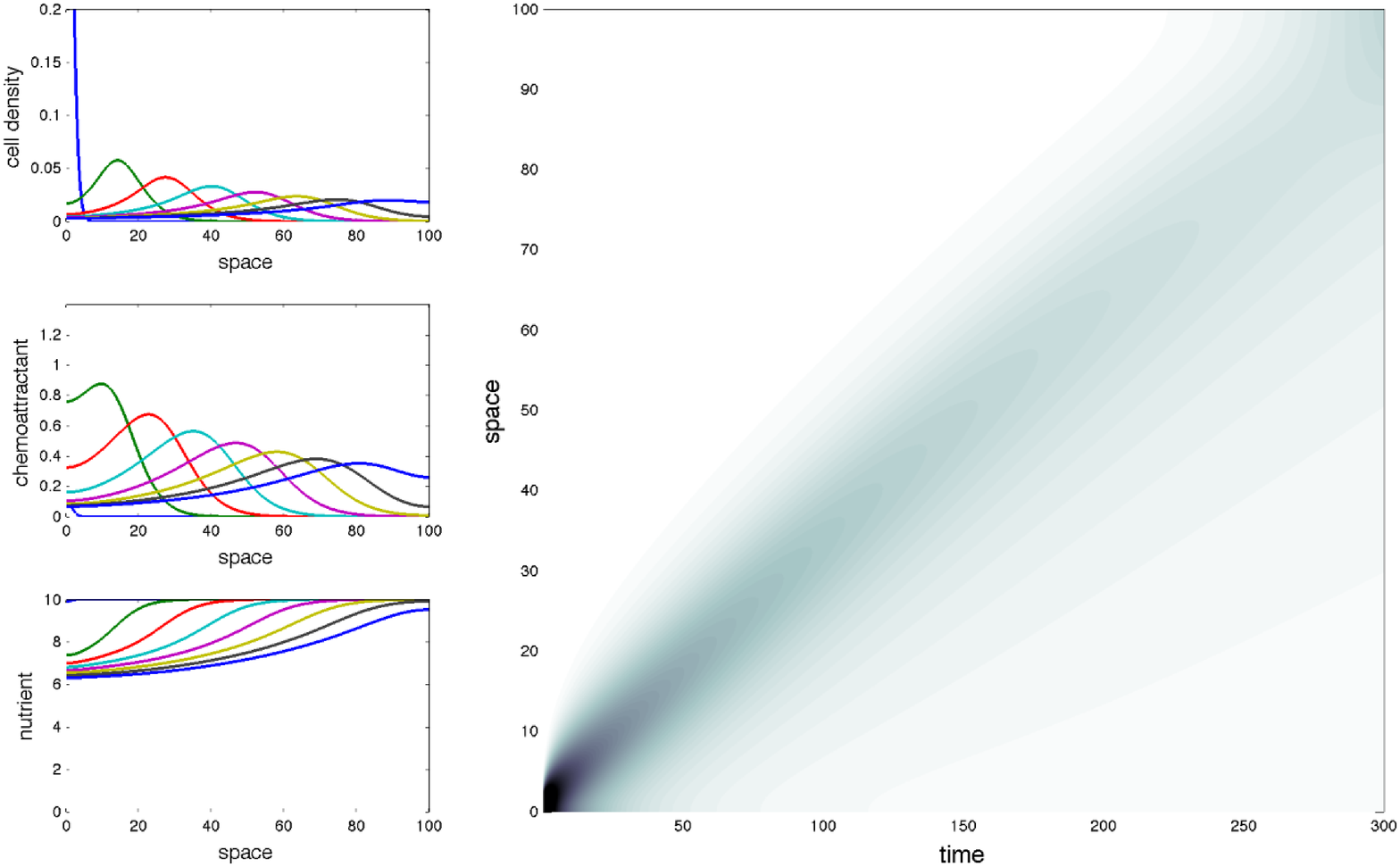}
\caption{\small 
(Top) In this experiment, bacteria are cultivated at a concentration of $5.10^8 \,{\rm cells/ml}$ in the same rich medium as in Fig. 3. After, they are resuspended in LB nutrient to an OD600 of $3.10^8\,{\rm cells/ml}$. We interpret the absence of pulse propagation as following. Bacteria are adapted to a rich environmnent before resuspension. Thus they are not able to sense small chemical fluctuations  necessary for clustering to occur when evolving in a relatively poor medium. 
(Bottom) Influence of the internal processes stiffness. When the individual response function $\phi_\delta$ is not stiff, the effect of dispersion is too strong and no pulse wave propagates, as opposed to  Fig.~\ref{fig:nolimited.sharp}.  Specific parameters are: $\delta = 10^{-1}$ and $N_0 = 10$ (see Section \ref{sec:num} for the other parameters).
In mathematical models of bacterial chemotaxis, it is commonly accepted that adaptation of cells to large chemoattractant changes acts through the measurement of relative time variations: $  S^{-1}{DS}/{Dt} $ (see Appendix). In our context, this is to say that the stiffness parameter $\delta$ is proportional to the chemical level $S$. Hence after having dramatically changed the environment and before bacteria adapt themselves, we can consider that the response function $\phi_\delta$ is not stiff. \label{fig:smooth}}
\end{center}
\end{figure}

\begin{figure}
\begin{center}
\includegraphics[width=\linewidth]{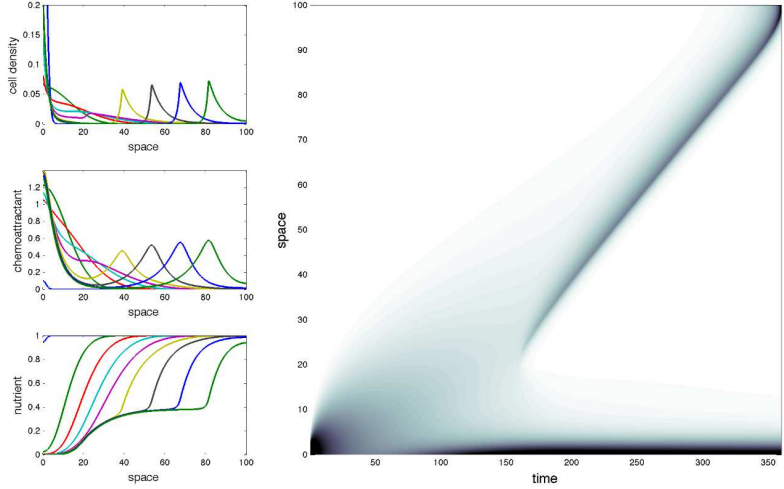} 
\caption{\small At low level of nutrient  the cell population splits into two subpopulations. A fraction remains trapped at the boundary (as a stationary profile) and a fraction  travels accross the channel with constant speed (see also Fig. \ref{fig:limited nutrient profile}).  Specific parameters are: $\delta = 10^{-3}$ and $N_0 = 1$ (see Section \ref{sec:num} for the other parameters).
}
\label{fig:limited nutrient}
\end{center}
\end{figure}

\begin{figure}
\begin{center}
\includegraphics[width=.48\linewidth]{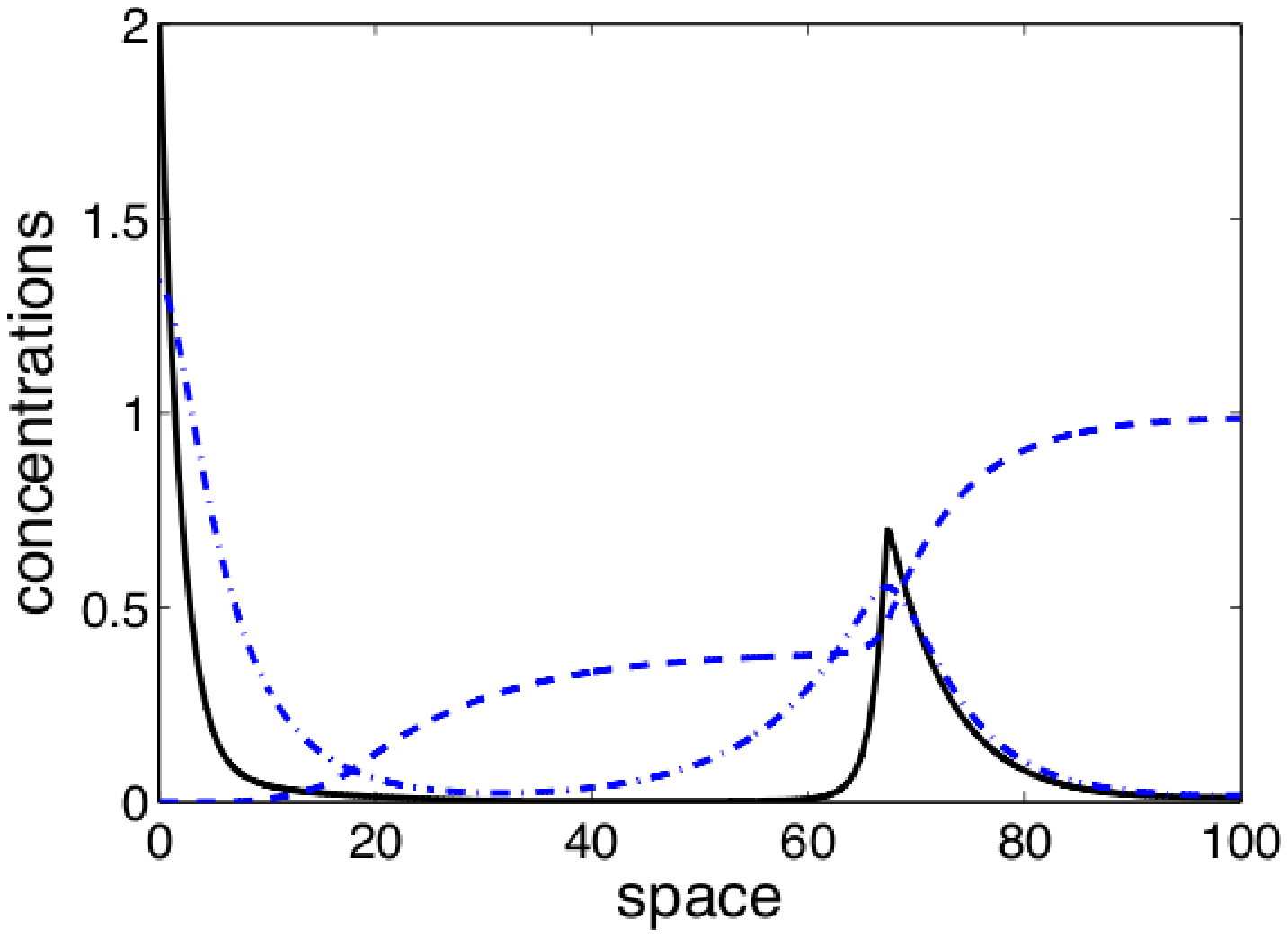} \, 
\includegraphics[width=.48\linewidth]{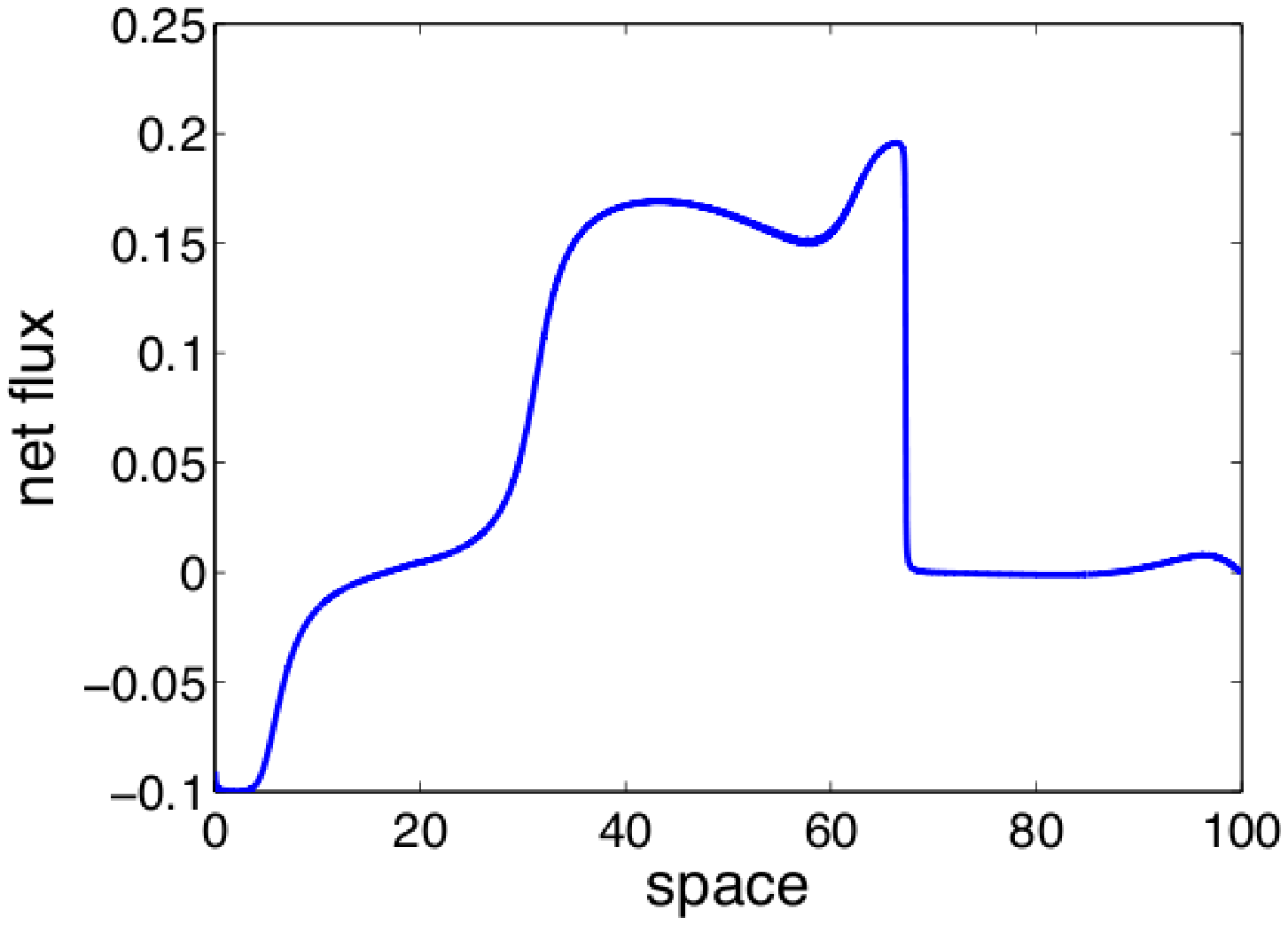} 
\caption{\small (Left) Respective profiles for the bacterial density (plain line), the chemical concentration (dash-dotted line) and the nutrient concentration (dashed line) in the limited nutrient/stiff response function configuration described in Fig. \ref{fig:limited nutrient} ($\delta = 10^{-3}$ and $N_0 = 1$). (Right) The net chemoattractive flux $u_S + u_N$ is plotted. It clearly shows a transition between the traveling wave (positive flux) and the equilibrium (negative flux) located at the left boundary.
}
\label{fig:limited nutrient profile}
\end{center}
\end{figure}

We opt for the following initial conditions in our numerical experiments: a decreasing exponential function centered on the left side of the channel for the cell density, no chemical signal, and a constant level of nutriment $N_0$.

Parameters are issued from literature (see \cite{Libchaber} and references therein) and from the mesoscopic derivation of the model (see Appendix 'Parameters and scales'). Adimensionalizing time and space yields the time and space scales being respectively  $10s$ and $200\mu m$, such that the total duration of the computation is approx. $1h$, and the length of the computational channel is $4cm$. 
We fix the following parameters, $D_\rho =1$, $\chi_S = \chi_N = 1$, for the cell density equation; $D_S = 2 $, $\alpha = 0.05$ for the chemical concentration \cite{Libchaber}. We assume in addition $D_N = 0$ for the sake of simplicity: the nutriment is not required to act as a non-local signal (and we expect this coefficient to have very little influence on the dynamics). By adimensionalizing further the system, we may choose $M = 1$, $\gamma = \beta = 1$ without loss of generality. Finally the signal response function $\phi_\delta$ is chosen as follows: $\phi_\delta(Y) = -  2/\pi \arctan(Y/\delta)$, with a stiffness factor $\delta$. We also keep the memory of the drift-diffusion limit performed in the Appendix, by setting $\eps = 0.1$. The only free parameters subject to variations are $N_0$ and $\delta$.

We can draw the following main conclusions from numerical experiments.



\subsubsection*{Influence of the stiffness of the internal response function (Fig. \ref{fig:smooth}).} When the stiffness assumption for the internal response function is relaxed, no pulse propagation is observed numerically. Deriving the exact conditions that guarantee the propagation of a traveling pulse seems to be a challenging task. However we give below some heuristics  for the particular choice $\phi_\delta(Y) = \phi(Y/\delta)$, $\delta$ being a stiffness factor.

Although the chemotactic equation of \eqref{eq:full model} is significantly different from the Keller-Segel model, they coincide as far as the stability of the homogeneous (unclustered) configuration is under question. We learn from Section \ref{sec:cluster} that the stiffness parameter $\delta$ plays an important role in the stability of the homogeneous solution. 
It is well known that the Keller-Segel system is subject to a bifurcation phenomenon due to its quadratic, non-local nonlinearity. This is well understood in two dimensions of space for instance \cite{PerthameBook,HillenPainter}. If some nondimensional parameter is small enough, diffusion dominates and no self-organization arises in the system. On the contrary, self-attraction between cells overcomes diffusion when this parameter is large, and yields the formation of a singularity ({\em i.e.} aggregation point) \cite{PerthameBook}.  

Clearly the same kind of mechanism acts here (see Fig. \ref{fig:nolimited.sharp} as opposed to Fig. \ref{fig:smooth}). However
there is no mathematical argumentation to sustain those numerical and intuitive evidence yet.  



\subsubsection*{Limited versus non limited nutrient (Fig. \ref{fig:limited nutrient}).}

When the nutrient is limited in the experimental device (and conditions for a pulse to travel are fulfilled) then only part of the bacterial population leaves the initial bump. The solution seems to be the superposition of a  traveling pulse and a stationary state (such as in Section \ref{sec:cluster} in the absence of nutrient). Solitary modes with smaller amplitudes may appear at the back of the leading one (not shown). To predict which fraction of mass starts traveling turns out to be a difficult question. 

\section{Conclusions and perspectives}

We present in this article a simple mathematical description for the collective motion of bacterial pulses with constant speed and asymmetric profile in a channel. The nature of this model significantly differs from the classical Keller-Segel system although it belongs to the same class of drift-diffusion equations. Our model is formally derived from a mesoscopic description of the bacterial density, which allows for a more accurate expression of the cell flux. We are able to compute analytically the speed of the pulse and its profile in the limit of a stiff response function $\phi_\delta$. The theoretical pulse speed has some striking features: it does not depend on the total number of bacteria, neither on the bacterial diffusion coefficient. This can be related to experimental evidence by Mittal {\em et al.} \cite{Mittal03} where bacteria self-organize into size-independent clusters. Our approach can be summarized as follows: a nutrient is added to pull chemotactic clusters of cells. This creates an imbalance in the fluxes which induces the asymmetry of the traveling profile. 


The next step would be to work at the kinetic level. Much has to be done for the design of efficient kinetic schemes for the collective motion of cells subject to chemotactic interactions. It would also be feasible to point out the dependency of the tumbling operator upon some internal variable ({\em e.g.} the cytoplasmic concentration of protein CheY). This approach carries out the coupling between an internal protein  network and the external chemoattractant signals \cite{SPO,Inoue07}.  Kinetic models are also relevant for describing this microscopic mechanism \cite{ErbanOthmer04,BournaveasC08} (the network is basically transported along the cells' trajectories). However the increase in complexity forces to reduce the size of the network, and to use rather caricatural systems mimicking high sensitivity to small temporal variations (excitation) and adaptation to constant levels of the chemoattractant.

Assuming independent integration of the chemical signals constitutes a strong hypothesis of our model. There exist two main membranous receptors triggering chemotaxis, namely Tar and Tsr. As the signals which act in the present experiments are not perfectly determined, 
we have considered the simplest configuration. To further analyse the interaction between the external signals, one should include more in-depth biological description of the competition for a single class of receptor \cite{Keymer}.

\appendix 
\section{Appendix: Kinetic models for chemotaxis and their drift-diffusion limit}
\label{sec:app}

\subsubsection*{Kinetic framework.}

The classical theory of drift-diffusion limit for kinetic modeling of bacterial chemotaxis is a way to compute the macroscopic fluxes $u_S$, $u_N$ in \eqref{eq:full model} \cite{HillenOthmer}. Because we assume a linear integration of the different signals for each individual, we restrict ourselves to the action of a single chemical species $S$.

The kinetic framework is as follows. A population of bacteria can be described at the mesoscopic scale by its local density $f(t,x,v)$ of cells located at the position $x$ and with velocity $v$ at time $t$. The kinetic equation proposed in the pioneering works of Alt, Dunbar and Othmer \cite{Ajmb80, ODA}  combines free runs at speed $v$, and tumbling events changing velocity from $v'$ (anterior) to $v$ (posterior), {\em resp.} according to the Boltzman type equation:
\begin{multline}\label{genkinmodel}
\displaystyle \partial_{t} f + v\cdot \nabla_{x} f = \int_{v'\in V} T[S](t,x,v'\to v) f(t,x,v')\, dv' \\ - \int_{v'\in V}T[S](t,x,v\to v') f(t,x,v)\, dv'   \ .
\end{multline}
The velocity space $V$ is  bounded and symmetric, usually $V = B(0,c)$ or  $V = \Sph(0,c)$ (bacteria having presumably constant speed). As we deal with the idealization of a two-dimensional phenomenon in one dimension of space, we shall perform our computations for $V = [-c,c]$, but the results contained in this paper do not depend on this particular choice. Kinetic models of chemotaxis have been studied 
recently in \cite{HKS,BCGP,BournaveasC08}.

\begin{figure}
\begin{center}
\includegraphics[width=.6\linewidth]{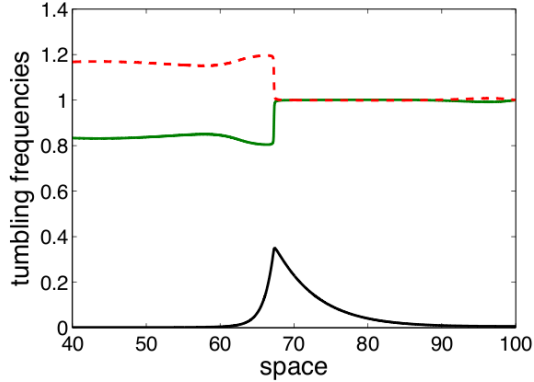} 
\caption{\small Tumbling frequencies (at the mesoscopic scale) obtained from the numerical experiment described in Fig. \ref{fig:limited nutrient}: the tumbling probability is higher when moving to the left (upper dashed line) at the back of the pulse, whereas the tumbling probability when moving to the right is lower (upper plain line), resulting in a net flux towards the right, as the pulse travels (see Fig. \ref{fig:limited nutrient profile}). Notice that these two curves are not symmetric w.r.t. to the basal rate 1, but the symmetry defect is of lower order ($10^{-3}$). The peak location is also shown for the sake of completeness (lower plain line). 
}
\label{fig:tumbling}
\end{center}
\end{figure}

The turning kernel $T$ describes the frequency of changing trajectories, from $v'$ to $v$. It expresses the way  external chemicals may influence  cell trajectories. A single bacterium is able to sense time variations of a chemical along its trajectory (through a time convolution whose kernel is well described since the experiments performed by Segall {\em et al.} \cite{Segall86}). For the sake of simplicity we neglect any memory effect, and we assume that a cell is able of sensing the variation of the chemical concentration along its trajectory. Following \cite{DolakSchmeiser}, this is to say that $T$ is given by the expression
\begin{equation} 
T[S](v'\to v) = \psi\left(\dfrac{DS}{Dt}\right) =  \psi\left(\partial_tS+v'\cdot \nabla_xS\right)\, . 
\label{eq:temporal response}
\end{equation}
The signal integration function $ \psi$ is non-negative and decreasing, expressing  that cells are less likely to tumble (thus perform longer runs) when the external chemical signal increases (see Fig. \ref{fig:tumbling} for such a tumbling kernel in the context of the present application). It is expected to have a stiff transition at 0, when the directional time derivative of the signal changes sign \cite{Segall86,SPO,Inoue07}.  Our study in Section \ref{sec:num}  boils down to  the influence of the stiffness, by introducing a one parameter family of functions $\psi_\delta(Y) = \psi(Y/\delta)$.

\subsubsection*{Parameters and scales.}

The main parameters of the model are the total number of bacteria $M$ which is conserved, the maximum speed of a single bacterium $c = \max\{ |v|; \; v\in V\}$, and the mean turning frequency $\lambda_0 = \psi_0 c^d$ (where $d$ denotes the dimension of space according to our discussion above). The main unknown is the speed of the traveling pulse, denoted by  $\sigma$.
We rescale the kinetic model \eqref{genkinmodel} into a nondimensional form as follows: 
\[ t = \tilde t \, \bar t \, , \quad  x = \tilde x\,  \bar x \, , \quad  v = \tilde v \,  c \, , \quad V = c \tilde V \, , \quad  T = \tilde T \,  \psi_0 \, . \]
We aim at describing traveling pulses in the regime $\bar x = \sigma \bar t$. Experimental evidence show that the bulk velocity $\sigma$ is much lower than the speed of a single bacterium $c$. This motivates to  introduce the ratio $\eps = \sigma/  c $. According to experimental measurements, we have $\eps \approx 0.1$. The kinetic equation writes: 
\begin{multline}
\label{eq:kinetic adim} \eps \partial_{\tilde t} \tilde f +  \tilde v\cdot \nabla_{\tilde x} \tilde f = \dfrac{\lambda_0 \bar x}{ c} \left\{\int_{\tilde v'\in \tilde V} \tilde \psi\left(\eps \partial_{\tilde t} \tilde S+  \tilde v'\cdot \nabla_{\tilde  x}\tilde S\right) \tilde f(\tilde t,\tilde x,\tilde v')\, d\tilde v' \right. \\ 
\left. - |\tilde V| \tilde \psi\left(\eps\partial_{\tilde t} \tilde S+  \tilde v\cdot \nabla_{\tilde x}\tilde S\right)  \tilde f(\tilde t,\tilde x,\tilde v)  \right\}  \ , 
 \end{multline}
 where $\tilde \psi (z) =  \psi(  c z/\bar x) $.
Following the experimental setting (see Introduction, Fig. \ref{fig:WaveChannel} and Fig. \ref{fig:front propagation}) and the biological knowledge \cite{BergBook}, we choose the scales $\bar x\approx 200\mu m$, $\lambda_0\approx 1s^{-1}$, and $ c = 20\mu m.s^{-1}$. Hence ${\lambda_0 \bar x}/{ c}\approx 10$. Therefore we rewrite this ratio as:
\[ 
\dfrac{\lambda_0 \bar x}{ c}  = \dfrac{\mu }{\eps}\, ,
\] 
where the nondimensional coefficient $\mu$ is of order 1. 

\subsubsection*{Drift-diffusion limit of kinetic models.}


To perform a drift-diffusion limit when $\eps\to 0$ (see \cite{HillenOthmer,CMPS,CDMOSS,PerthameBook}, and \cite{DolakSchmeiser,FilbetLaurencotPerthame} for other scaling limits, {\em e.g.} hyperbolic), we shall assume that the variations of $\psi$ around its meanvalue $\psi_0$ are of amplitude $\eps$ at most. It writes in the nondimensional version as follows: $\psi(Y) = 1 + \eps \phi_\delta(Y)$. Hence the chemotactic contribution is a perturbation of order $\eps$ of a unbiased process which is constant in our case because the turning kernel does not depend on the posterior velocity and the first order contribution is required to be symmetric with respect to $(v',v)$. 
This hypothesis is in agreement with early biological measurements. 
It is also relevant from the mathematical viewpoint as we are looking for a traveling pulse regime where the speed of the expected pulse is much slower than the speed of a single individual. This argues in favour of a parabolic scaling as performed in this Appendix. 





The rest of this Appendix is devoted to the derivation of the Keller-Segel type model in one dimension of space:
\begin{equation} \label{eq:parabolic}
\displaystyle \partial_t \rho + \partial_{x}\left( - D_\rho \partial_x \rho + \rho u_S \right) = 0 \,  .
\end{equation}
Dislike the classical Keller-Segel model (used for instance by Salman et al. \cite{Libchaber}), singularities  cannot form (excessively populated aggregates) with the  chemotactic flux $u_S$ given in \eqref{eq:kinflux} below. This is because the latter remains uniformly bounded (see also Mittal {\em et al.} \cite{Mittal03} where clusters emerge which are plateaux and thus not as singular as described for KS system in a mathematical sense). 


\begin{proof}[Sketch of parabolic derivation.]
We start from the nondimensional kinetic equation \eqref{eq:kinetic adim}:
\begin{multline*} 
\eps \partial_t f +  v\cdot \nabla_x f = \dfrac{\mu}\eps\left\{\int_{v'\in V} \left(1 + \eps \phi_\delta[S](v')\right) f(t,x,v')\, dv' \right.\\ 
\left. - |V|\left( 1+ \eps \phi_\delta[S](v)\right) f(t,x,v)\right\}\, , 
\end{multline*}
which reads as follows,
\begin{multline}
\eps \partial_t f +  v\cdot \nabla_x f =  \dfrac{\mu}\eps \left(\rho(t,x) - |V|   f(t,x,v)\right) \\ +  \mu\left( \int_{v'\in V} \phi_\delta[S](v') f(t,x,v')\, dv' - |V| \phi_\delta[S](v)  f(t,x,v)\right)\, .
\label{eq:parabolic scaling}
\end{multline}
Therefore the dominant contribution in the tumbling operator is a relaxation towards a uniform distribution in velocity at each position: $f(t,x,v) = \rho(t,x) F(v)$ as $\eps\to 0$, where $F(v) = |V|^{-1}\ind{v\in V}$. Notice that more involved velocity profiles can be handled \cite{CMPS,Perthame04}, but this is irrelevant in our setting as the tumbling frequency does not depend on the posterior velocity $v$.
 
The space density $\rho(t,x)$ remains to be determined. For this purpose we first integrate with respect to velocity $v$ and we obtain the equation of motion for the local density $\rho(t,x) = \int_{v\in V}f(t,x,v)dv$:
\[\partial_t \rho + \nabla\cdot j  = 0\, ,\quad j = \eps^{-1} \int_{v\in V} v f(t,x,v)\, dv\, .\]
To determine the bacterial flow $j$ we integrate \eqref{eq:parabolic scaling} against $v$:
\begin{multline*}  \eps \partial_t\left( \int_{v\in V} v f(t,x,v)\, dv \right)  + \nabla_x\cdot\left( \int_{v\in V} v\otimes v f(t,x,v)\, dv \right)\\ =  - \mu |V| j  - \mu|V| \int_{v\in V} v \phi_\delta[S](v) f(t,x,v)\, dv\, .
\end{multline*}
We obtain formally, as $\eps \to 0$: 
\begin{equation}
j = -   \nabla_x  \left(  \rho(t,x)  \dfrac1{\mu d |V|^2} \int_{v\in V} |v|^2 \, dv \right) - \rho(t,x)  \dfrac1{|V|} \int_{v\in V} v \phi_\delta[S](v) \, dv  \,  . 
\end{equation}
Finally, the drift-diffusion limit equation \eqref{eq:parabolic} reads in one dimension of space:
\begin{equation}
\partial_t \rho = \left(\dfrac{1}{4\mu}\int_{v\in [-1,1]} |v|^2\, dv\right) \partial^2_{xx} \rho +  \partial_x  \left(\rho \int_{v\in [-1,1]} v \phi_\delta\left( \eps\partial_t S + v\partial_x S\right) \, \frac{dv}{2}  \right)\, .
\end{equation}
\end{proof}

To sum up, we have derived a macroscopic drift-diffusion equation, where the bacterial diffusion coefficient and the chemotactic flux are given by:
\begin{equation} D_\rho = \dfrac{1}{4\mu}\int_{v\in [-1,1]} |v|^2\, dv\, , \quad u_S = - \int_{v\in [-1,1]} v \phi_\delta\left( \eps\partial_t S + v\partial_x S\right) \, \frac{dv}{2} \, . 
\label{eq:kinflux}
\end{equation}
In the limiting case  where the internal response function  $\phi_\delta$ is bivaluated: $\phi_\delta(Y) =  \phi_0 \ind{Y<0} -  \phi_0 \ind{Y>0}$, the flux rewrites simply as: 
\[u_S = \dfrac{\phi_0}2  \left(1 -  \left(\epsilon\dfrac{ \partial_t S}{ \partial_x S}\right)^2  \right)_+ \sign(\partial_x S)\, . 
\]

For the sake of comparison, we highlight the corresponding expressions which have been obtained by Dolak and Schmeiser. In \cite{DolakSchmeiser} authors perform a hyperbolic scaling limit leading to the following chemotactic equation for the density of bacteria:
\[  
\partial_t \rho + \nabla\cdot\left( - \eps D\nabla \rho  + \rho   U_S \right) = 0 \,  ,
 \]
where $D$ is an anisotropic diffusion tensor and the   chemotactic flux is given by:
\[   U_S = \left(\dfrac1A \int_{v\in V} \dfrac {v_1}{\psi(\partial_t S + v_1\nabla S|)}\, dv\right) \dfrac{\nabla S}{|\nabla S|} \,  , \]
for some renormalizing factor $A$.
The two approaches do not differ that much at first glance (especially when $\psi$ is bivaluated). Notice however that the "small" $\epsilon$ parameter does not appear at the same location: in front of the diffusion coefficient in the hyperbolic limit and inside the chemotactic flux in the parabolic limit.

\end{document}